# Notes on the occupancy problem with infinitely many boxes: general asymptotics and power laws*

**Alexander Gnedin**

*Utrecht University*
*e-mail:* gnedin@math.uu.nl

**Ben Hansen**

*University of Michigan*
*e-mail:* bbh@umich.edu

**Jim Pitman**

*University of California, Berkeley*
*e-mail:* pitman@stat.Berkeley.EDU

**Abstract:** This paper collects facts about the number of occupied boxes in the classical balls-in-boxes occupancy scheme with infinitely many positive frequencies: equivalently, about the number of species represented in samples from populations with infinitely many species. We present moments of this random variable, discuss asymptotic relations among them and with related random variables, and draw connections with regular variation, which appears in various manifestations.

## Contents

*Research supported in part by N.S.F. Grant DMS-0405779 and N.I.C.H.D. Grant HD045753-01

## 1. Introduction

We consider the classical multinomial occupancy scheme in which balls are thrown independently at a fixed infinite series of boxes, with probability $p_j$ of hitting the $j$th box. The frequencies $(p_j, \ j = 1, 2, \ldots)$ are assumed nonincreasing, strictly positive and satisfying $\sum_j p_j = 1$. As $n$ balls are thrown, their allocation is captured by the array $X_n = (X_{n,j}, \ j = 1, 2, \ldots)$, where $X_{n,j}$ is the number of balls out of the first $n$ that fall in box $j$.

In concrete applications, instead of boxes one has types or species of sampling units, and the sample array of types, $X_n$, is of interest for what it reveals about the population frequencies $(p_j)$. Such *species sampling problems* arise in ecology, to be sure, but also in database query optimization, where the sampling units may be entries in columns of a database while the species consist of all of distinct values appearing in the column [10]; in literature, where the sampling units may be words appearing in a given author's known works while the species consist of all words known to that author [14]; in disclosure risk limitation, where the sampling units may be people or firms listed in a microdata file, without names or other overtly identifying information, while the types are unique combinations of values of variables with which the people or firms might be implicitly identified [34]; and in many other areas [9]. Models positing infinitely many boxes or species may approximate sampling from large, finite populations, or they may be useful as models of 'superpopulations' from which both samples and background populations are notionally drawn.

A functional of $X_n$ which appears in many contexts is the number of nonempty boxes
$$K_n = \#\{j: \ X_{n,j} > 0\}.$$
$K_n$ is sometimes regarded as a measure of diversity of the sample. More detailed information is carried by the counts of boxes occupied by exactly $r$ balls
$$K_{n,r} = \#\{j: \ X_{n,j} = r\} \quad (r = 1, 2, \ldots),$$
so that $K_n = \sum_r K_{n,r}$ and $\sum_r r\, K_{n,r} = n$. The combinatorial object encoded into the array of counts $(K_{n,1}, \ldots, K_{n,n})$ is a random partition of integer $n$; this partition has $K_n$ parts which correspond to positive entries of $X_n$.

The variables $K_n$ and $K_{n,r}$'s have been intensively studied in the occupancy scheme with finitely many positive frequencies, which in some models may vary in a certain way with $n$. The most studied is, of course, the classical case of $m$ equal frequencies (e.g. in the familiar 'birthday paradox' one is interested in the probability of the event $K_{n,1} < n$). Kolchin, Sevast'yanov and Chistyakov [29] identified five distinct asymptotic regimes for $m = m(n) \to \infty$ to secure either a Poisson or a normal limit distribution for $K_n$ (more precisely, they discussed the number of empty boxes $m - K_n$); for example, when $m = n$, both the mean and the variance of $K_n$ grow approximately linearly with $n$, and the limit distribution is normal. See [27, 26, 29] for extensions, surveys and many references.

In contrast to that, the literature on the problem with infinitely many fixed positive frequencies is rare. For fixed frequencies $(p_j)$ the asymptotic growth of each $X_{nj}$ (as $n \to \infty$) is linear, and the growth of $K_n$ is always sublinear since the main contribution to $K_n$ occurs due to the boxes whose frequencies are arbitrarily small. More delicate asymptotics of the moments of $K_n$ and $K_{n,r}$ are determined by the way $p_j$ approach 0. The first systematic study of the asymptotics appeared in a remarkable paper by Karlin [28], in which he proved a central limit theorem under a condition of regular variation on the frequencies $(p_j)$.

Recently, two new sources of interest to the infinite model emerged. On the one hand, it has been observed that 'power laws' for $K_n$ are quite common in partition-valued processes of coagulation and fragmentation (see [32, 22, 24, 5]); these fall in the range of regular variation with positive index. On the other hand, the case of geometric sequence $(p_j)$ (which is an instance of slow variation) was intensively studied in connection with the analysis of algorithms [2, 30].

These notes present a kind of a survey which extends and updates Karlin's results on the infinite occupancy scheme. The results in Sections 2-6 and 10 are of general nature, while in Sections 7-9 we work under the assumption of regular variation. In particular, we record various guises of regular variation, and mention some recent developments. Some of the results are new and others are scattered in the literature and have been several times re-discovered (especially the results in Section 10).

*Notation.* Throughout $c, c_j$ denote positive constants whose values are not important and may change from line to line. We use $f \sim g$, $f \ll g$, $f \gg g$ and $f \asymp g$ for $f/g \to 1$, $f/g \to 0$, $f/g \to \infty$ and $c_1 < f/g < c_2$, respectively. When one either of $f$ and $g$ is a random quantity, the notation $f \sim_{\text{a.s}} g$, $f \ll_{\text{a.s}} g$ etc. means that the asymptotic relation holds with probability one. Convergence in distribution is denoted $\to_d$.

## 2. Moments and poissonization

We start by recalling the familiar fact that $X_n := (X_{n,1}, X_{n,2}, \ldots)$ has a multinomial distribution with parameters $(n, (p_j))$:

$$\mathbb{P}(X_{n,j} = n_j, \ j = 1, 2, \ldots) = \frac{n!}{\prod_j n_j!} \prod_j p_j^{n_j} \qquad (n = \sum_j n_j).$$

From this the distribution of partition $(K_{n,1}, \ldots, K_{n,n})$ is recovered by summation. Specifically, for $(k_1, \ldots, k_n)$ a fixed partition of $n$

$$\mathbb{P}((K_{n,1}, \ldots, K_{n,n}) = (k_1, \ldots, k_n)) = \frac{n!}{n_1! \cdots n_k!} \sum_{\substack{\text{distinct} \\ j_1, \ldots, j_k}} p_{j_1}^{n_1} \cdots p_{j_k}^{n_k},$$

where $n_1, \ldots, n_k$ is a sequence of length $k = \sum_i k_i$, with $k_i$ terms equal to $i$ for $i = 1, \ldots, n$. The infinite sum is called the monomial symmetric function in the variables $p_j$. Formulas for the distribution of $K_n$ and marginal distributions of $K_{n,r}$'s follow by further summation over partitions of $n$. In terms of the generating function, the probability of $K_n = m$ is equal to the coefficient at $x^n y^m / n!$ in the series expansion of the infinite product

$$\prod_j (1 + y(e^{p_j x} - 1)).$$

Formulas for the moments follow from the representation

$$K_{n,r} = \sum_j \mathbf{1}(X_{n,j} = r), \quad K_n = \sum_j \mathbf{1}(X_{n,j} > 0), \tag{1}$$

(where $\mathbf{1}(\cdots)$ equals 1 if $\cdots$ is true and equals 0 otherwise). Denoting

$$\Phi_n := \mathbb{E}[K_n], \quad \Phi_{n,r} := \mathbb{E}[K_{n,r}]$$

we easily see that

$$\Phi_n = \sum_j (1 - (1-p_j)^n), \quad \Phi_{n,r} = \binom{n}{r} \sum_j p_j^r (1-p_j)^{n-r}. \tag{2}$$

These are related by the formulas

$$\Phi_n = \sum_r \Phi_{n,r}, \quad \Phi_{n,r} = (-1)^{r-1} \binom{n}{r} \Delta^r \Phi_n,$$

where $\Delta^r$ is the $r$th iterate of the difference operator $\Delta \Phi_n = \Phi_n - \Phi_{n-1}$. Lengthy but straightforward computations yield formulas for the variance

$$V_n := \mathrm{Var}[K_n] = \Phi_{2n} - \Phi_n + \sum_{i \neq j} [(1 - p_i - p_j)^n - (1-p_i)^n (1-p_j)^n], \tag{3}$$

$$V_{n,r} := \mathrm{Var}[K_{n,r}] = \Phi_{n,r} - \frac{\binom{n}{r}\binom{n}{r}}{\binom{2n}{2r}} \Phi(2n, 2r) + \tag{4}$$

$$\binom{n}{r} \sum_{i \neq j} p_i^r p_j^r \left[ \binom{n-r}{r} (1 - p_i - p_j)^{n-2r} - \binom{n}{r} (1-p_i)^{n-r} (1-p_j)^{n-r} \right].$$

Formulas for higher moments and covariance can be derived in the same way, but they seem to be of little practical use in the explicit form. See [29, Section 3.1] for such moments computations, which are also valid in the case of infinitely many positive frequencies.

One major obstacle in the study of counts $K_n$ and $K_{n,r}$'s is that the indicators in (1) are not independent. A common recipe to circumvent this difficulty is to first consider a closely related type of model in which the balls are thrown in continuous time at epochs of a unit rate Poisson process $(P(t), t \geq 0)$, which is independent of $(X_n, n = 1, 2, \ldots)$. The advantage of this randomization is that one can exploit independence, since the balls fall then in the boxes according to independent Poisson processes $(X_j(t), t \geq 0)$, at rate $p_j$ for box $j$. Once the properties of the Poisson allocation scheme are acquired, one still needs to translate them in the fixed-$n$ results, by using a kind of depoissonization technique.

**Remark.** Poissonization liberates us from the constraints of the fixed-$n$ scheme. The Poisson model is well defined for arbitrary positive rates $p_j$ which need not satisfy $\sum_j p_j < \infty$ or even $p_j < 1$. If $\sum_j p_j < \infty$, a reduction to the normalized case $\sum_j p_j = 1$ is maintained by the obvious time-change $t \mapsto t/\sum_j p_j$. If $\sum_j p_j = \infty$, $K(t)$ is infinite with probability one, though $K_r(t)$'s can be still finite. One can also consider two-sided infinite sequences like geometric ($p_j = q^j$, $j \in \mathbb{Z}$) (with $0 < q < 1$), for which the $K_r(t)$'s and the sum of variances $V(t) := \sum_{j=-\infty}^{\infty} \text{Var}\,[\mathbf{1}(X(t) > 0)] = \sum_{j=-\infty}^{\infty} (e^{-p_j t} - e^{-2p_j t})$ are all finite.

The convention in this paper is that the quantities associated with the Poisson allocation scheme appear in the functional notation, while the lower-index notation is reserved for the fixed-$n$ scheme. For instance, $X_j(t) = X_{P(t),j}$,

$$K(t) := K_{P(t)} = \sum_j \mathbf{1}(X_j(t) > 0), \quad K_r(t) := K_{P(t),r} = \sum_j \mathbf{1}(X_j(t) = r). \quad (5)$$

For the poissonized moments $\Phi(t) := \mathbb{E}\,[K(t)]$, $\Phi_r(t) := \mathbb{E}\,[K_r(t)]$ we have

$$\Phi(t) = \sum_j (1 - e^{-tp_j}), \quad \Phi_r(t) = \frac{t^r}{r!} \sum_j p_j^r \, e^{-tp_j} \qquad (r = 1, 2, \ldots),$$

and these are related via

$$\Phi(t) = \sum_r \Phi_r(t), \quad \Phi_r(t) = \frac{t^r}{r!} \Phi^{(r)}(t),$$

where $\Phi^{(r)}(t)$ is the $r$th derivative. Formulas for the variance are simpler than (3) and (4) because the cross-terms disappear due to independence:

$$V(t) := \text{Var}\,[K(t)] = \Phi(2t) - \Phi(t), \quad V_r(t) := \text{Var}\,[K_r(t)] = \Phi_r(t) - 2^{-2r}\binom{2r}{r}\Phi_{2r}(2t). \tag{6}$$

Similarly, for integer $r \neq s$, using $K_r(t) + K_s(t) = \sum_j \mathbf{1}(X_j(t) \in \{r, s\})$ the covariance is computed as

$$\operatorname{cov}(K_r(t), K_s(t)) = -2^{-r-s}\binom{r+s}{r}\Phi_{r+s}(2t). \tag{7}$$

For analytical reasons which will be soon clear it is convenient to encode the frequencies $(p_j)$ into an infinite counting measure

$$\nu(\mathrm{d}x) := \sum_j \delta_{p_j}(\mathrm{d}x)$$

on $]0, 1[$, where $\delta_x$ is the Dirac mass at $x$. Equivalently,

$$\sum_j f(p_j) = \int_0^1 f(x)\nu(\mathrm{d}x) \tag{8}$$

holds for arbitrary $f \geq 0$. In particular,

$$\Phi_n = \int_0^1 (1 - (1-x)^n)\nu(\mathrm{d}x), \tag{9}$$

$$\Phi_{n,r} = \binom{n}{r}\int_0^1 x^r(1-x)^{n-r}\nu(\mathrm{d}x) \quad (r = 1, 2, \ldots), \tag{10}$$

$$\Phi(t) = \int_0^1 (1 - e^{-tx})\nu(\mathrm{d}x), \tag{11}$$

$$\Phi_r(t) = \frac{t^r}{r!}\int_0^1 x^r e^{-tx}\nu(\mathrm{d}x) \quad (r = 1, 2, \ldots). \tag{12}$$

**Remarks.** The formulas for expected values remain exactly the same when the frequencies $(p_j)$ are random, in which case the 'intensity measure' $\nu$ is defined by taking expectation in the left-hand side of (8). See the recent work on composition structures [22, 23, 4] for more in this direction.

The summability constraint $\sum_j p_j = 1$ translates as $\int_0^1 x\,\nu(\mathrm{d}x) = 1$, and implies that $\nu$ can be also interpreted as a Lévy measure of some subordinator, which jumps by $p_j$ at rate 1 for each $j$. In this interpretation (11) has the meaning of a Laplace exponent [20].

Both $\Phi(t)$ and $\Phi_n$ (considered as functions of a real or complex-valued argument) are Bernstein functions which uniquely determine $\nu$, see [20]. They are related by the poissonization identity

$$\Phi(t) = \mathbb{E}[\Phi_{P(t)}] = \sum_{n=0}^{\infty} \frac{e^{-t}t^n}{n!}\Phi_n \quad (t \geq 0).$$

## 3. Some estimates of the moments

Monotonicity is a key feature of $K_n$. In fact, we have $K_n \uparrow_{\text{a.s.}} \infty$ (as $n \to \infty$) and $K(t) \uparrow_{\text{a.s.}} \infty$ (as $t \to \infty$), because each box is eventually discovered by a ball.

By monotone convergence, also $\Phi_n \uparrow \infty$ and $\Phi(t) \uparrow \infty$. However, the growth is always sublinear: $\Phi_n \ll n$ ($n \to \infty$) and $\Phi(t) \ll t$ ($t \to \infty$). Indeed, if we ignore the first $J$ boxes, the mean number of discovered boxes among the remaining ones is at most $n \sum_{j>J} p_j$ (respectively, at most $t \sum_{j>J} p_j$ in the Poisson scheme). Thus $\Phi_n < J + n \sum_{j>J} p_j$ (respectively, $\Phi(t) < J + t \sum_{j>J} p_j$), and selecting $J$ arbitrarily large we see that $\limsup \Phi(n)/n = \limsup \Phi(t)/t = 0$. Aside from these two general features, the growth properties of $K_n$ can be fairly arbitrary.

The next lemma gives general estimates of closeness of the moments in the fixed-$n$ scheme and the Poisson scheme.

**Lemma 1.** *For $n \to \infty$ the following estimates hold:*

$$
\begin{aligned}
|\Phi(n) - \Phi_n| &< \frac{2}{n} \Phi_2(n) \to 0, \\
|\Phi_r(n) - \Phi_{n,r}| &< \frac{c}{n} \max\{\Phi_r(n), \Phi_{r+2}(n)\} \to 0, \\
|V(n) - V_n| &< \frac{c}{n} \max\{\Phi_1(n), \Phi_1(n)^2\} < \frac{c}{n} \max\{1, \Phi_1(n)^2\}, \\
|V_r(n) - V_{n,r}| &< \frac{c}{n} \max\{\Phi_r(n)^2, \Phi_{r+1}(n)^2, \Phi_{r+1}(n), \Phi_{r+2}(n), \Phi_{2r}(2n)\}.
\end{aligned}
$$

*Proof.* The first two bounds follow from the elementary inequality $0 \leq e^{-nx} - (1-x)^n \leq nx^2 e^{-nx}$ valid for $0 \leq x \leq 1$. For instance,

$$
\begin{aligned}
|\Phi(n) - \Phi_n| &= \int_0^1 \left(e^{-nx} - (1-x)^n\right) \nu(\mathrm{d}x) < n \int_0^1 x^2 e^{-nx} \nu(\mathrm{d}x) \\
&= \frac{2}{n} \Phi_2(n) < \frac{2}{n} \Phi(n) \to 0.
\end{aligned}
$$

The last two bounds follow from this and estimates of the cross-terms in (3) and (4), by using the expansion

$$(a-b)^n = a^n + n a^{n-1} b + O(n^2 a^{n-2} b^2) \quad (a > b)$$

with $a = (1-p_i)(1-p_j) = (1 - p_i - p_j + p_i p_j)$ and $b = p_i p_j$. $\square$

Taken together with $\Phi(t) \uparrow \infty$ the lemma implies $\Phi_n \sim \Phi(n)$.

**Remark.** It seems plausible that the relations $V(n) \sim V_n$ and $\Phi_r(n) \sim \Phi_{n,r}$ are true for arbitrary $(p_j)$. However, the estimates in the lemma are not strong enough to entail such a conclusion in full generality, although it has been shown under various circumstances (see e.g. Lemma 4 below and Section 6) and no counterexamples are known. Hwang and Janson [25, Proposition 4.3(ii)] show that always $V(n) \asymp V_n$. The difficulty is that $V(t)$ and $\Phi_r(t)$ may exhibit rather irregular oscillatory behaviour. For instance $V(t)$, $V_n$, $\Phi_r(t)$ and $\Phi_{n,r}$ may approach 0 arbitrarily closely for some $n$ or $t$, (though they cannot converge to 0, as is seen by selecting a subsequence $n_j \asymp 1/p_j$ or $t_j \asymp 1/p_j$, respectively).

We denote
$$\vec{\nu}(x) := \nu[x, \infty[$$
the right tail of $\nu$. Note that there are at most $m$ frequencies not smaller than $1/m$, hence $\vec{\nu}(1/m) < m$. Moreover, we have $\vec{\nu}(x) \ll x^{-1}$ for $x \downarrow 0$. Indeed, integrating by parts for $x > 0$ we obtain

$$\int_x^1 u\nu(\mathrm{d}u) = x\,\vec{\nu}(x) + \int_x^1 \vec{\nu}(u)\mathrm{d}u.$$

As $x \downarrow 0$ the integral at left increases to 1, entailing by monotonicity that the integral at right converges, and by extension that $x\,\vec{\nu}(x)$ converges to a limit also. Since $\lim x\,\vec{\nu}(x) > 0$ would force the integral at right to diverge, $x\,\vec{\nu}(x) \to 0$.

## 4. Laws of large numbers

The mean number of occupied boxes satisfies $\Phi(t+\tau)-\Phi(t) < \Phi(\tau)$ (for $\tau, t > 0$). One way to justify this is by noting that the mean number of distinct boxes hit during any time interval $[\tau, \tau+t]$ is $\Phi(t)$, but some of them have been discovered before time $t$ and do not contribute to $K(t+\tau)$. The same follows from concavity of $\Phi(t)$, which implies

$$\frac{\Phi(t+\tau) - \Phi(t)}{\tau} < \frac{\Phi(\tau) - \Phi(0)}{\tau}.$$

Similar inequalities hold for $\Phi_n$. Using these the variance can be bounded via expectation as

$$V(t) = \Phi(2t) - \Phi(t) < \Phi(t) = \mathbb{E}\left[K(t)\right], \quad V_n < \Phi_{2n} - \Phi_n < \Phi_n = \mathbb{E}\left[K_n\right].$$

Applying Chebyshev's inequality and recalling that $\Phi_n \uparrow \infty$ and $\Phi(t) \uparrow \infty$, the bound on the variance allows one to conclude that both $K(t)/\Phi(t)$ and $K_n/\Phi_n$ converge to 1 in probability, which is a result due to Bahadur [3]. A similar analysis invoking (6), (4) and Lemma 1 shows that also $K_r(t)/\Phi_r(t)$ and $K_{n,r}/\Phi_{n,r}$ converge to 1 in probability, provided $\Phi_r(t) \to \infty$. For $K_n$ and $K(t)$ there is the following strengthening due to Karlin [28, Theorem 8].

**Proposition 2.** *For arbitrary* $(p_j)$

$$K_n \sim_{\text{a.s.}} \Phi_n, \quad K(t) \sim_{\text{a.s.}} \Phi(t).$$

*Proof.* The function $\Phi(t)$ is continuous, increasing and satisfies $\Phi'(t) < 1$. Thus it is possible to select an increasing sequence $(t_m, \ m = 1, 2, \ldots)$ such that $m^2 < \Phi(t_m) < m^2 + 1$. We have then from the above estimate of the variance $\mathbb{P}(|K(t_m)/\Phi(t_m) - 1| > \epsilon) < \epsilon^{-2} m^{-2}$, and by summability of the bound $K(t_m)/\Phi(t_m) \to_{\text{a.s.}} 1$ along the subsequence. For $t_m < t < t_{m+1}$, the monotonicity implies $K(t_m) \leq K(t) \leq K(t_{m+1})$ and $\Phi(t_m) < \Phi(t) < \Phi(t_{m+1})$. This allows one to squeeze the ratio as

$$\frac{K(t_m)}{\Phi(t_{m+1})} < \frac{K(t)}{\Phi(t)} < \frac{K(t_{m+1})}{\Phi(t_m)},$$

where both sides converge to 1 almost surely, in consequence of the above and $\Phi(t_m)/\Phi(t_{m+1}) \to 1$. The argument for $K_n$ is completely analogous. □

Instead of using the Chebyshev inequality one can exploit a finer Bernstein-type large deviation bound for sums of independent bounded variables [16, p. 911]:

$$\mathbb{P}(|K(t)/\Phi(t) - 1| \geq \epsilon) < c\, e^{-\epsilon' \Phi(t)},$$

where $\epsilon'$ depends on $\epsilon$. This allows one to choose a subsequence $\{t_m\}$ with smaller gaps, so that $\Phi(t_{m+1}) - \Phi(t_m) \asymp m^{-1}$.

Using monotonicity of $\sum_{r \geq s} K_{n,r}$, we obtain along the same lines for every fixed integer $s$

$$\sum_{r \geq s} K_{n,r} \sim_{\text{a.s.}} \sum_{r \geq s} \Phi_{n,r}, \quad \sum_{r \geq s} K_r(t) \sim_{\text{a.s.}} \sum_{r \geq s} \Phi_r(t) \ \text{ and } \ \sum_{r \geq s} \Phi_{n,r} \sim \sum_{r \geq s} \Phi_r(t).$$

**Remark.** The relations $K_{n,r} \sim_{\text{a.s.}} \Phi_{n,r}, K_r(t) \sim_{\text{a.s.}} \Phi_r(t)$ may fail, simply because $\Phi_r(t)$ need not go to $\infty$, while the counting processes have unit jumps. It is natural to conjecture that these laws of large numbers are true under the condition $\Phi_r(t) \to \infty$, but we do not know if this has been proved in full generality.

## 5. CLT for $K_n$

Recall the representation (5) of $K(t)$ as a sum of independent indicators. Applying the Lindeberg-Feller condition, we see that $(K(t) - \Phi(t))/V(t)^{1/2}$ converges to the standard normal distribution provided $V(t) \to \infty$. The following depoissonization argument leading to the CLT for $K_n$ follows the line in Dutko [12]. Monotonicity of $K_n$ plays here a central role.

**Proposition 3.** *The conditions $V(t) \to \infty$ and $\text{Var}[K_n] \to \infty$ are equivalent, and if they hold, the law of $(K_n - a_n)/b_n^{1/2}$ converges to the standard normal distribution, where $\Phi_n$ or $\Phi(n)$ can be selected for the constant $a_n$ and $V(n)$ or $V_n$ for $b_n$.*

The next lemma will imply that both choices for $b_n$ are good.

**Lemma 4.** *The conditions* $\operatorname{Var}[K_n] \to \infty$ *and* $V(t) \to \infty$ *are equivalent, and imply* $V(n) \sim V_n$.

*Proof.* We first show that $\Phi_1(t)^2 \ll tV(t)$. Denote for shorthand

$$\phi(t) := t^{-1}\Phi_1(t) = \Phi'(t).$$

Note that $\phi'(t) < 0$. Since $\vec{\nu}(x) < \infty$ for $x > 0$ we have

$$\phi(t) \sim \int_0^\epsilon xe^{-tx}\nu(\mathrm{d}x) \quad (t \to \infty)$$

for every $\epsilon > 0$. Thus by Cauchy-Schwarz (applied to the measure $x\nu(\mathrm{d}x)$)

$$\phi(t)^2 \sim \left[\int_0^\epsilon xe^{-tx}\nu(\mathrm{d}x)\right]^2 \leq \int_0^\epsilon x\nu(\mathrm{d}x) \int_0^\epsilon xe^{-2tx}\nu(\mathrm{d}x) \sim \phi(2t) \int_0^\epsilon x\nu(\mathrm{d}x),$$

and letting $\epsilon \to 0$ the first integral factor vanishes, which yields $\phi(t)^2 \ll \phi(2t)$. Using that $\phi$ is decreasing,

$$V(t) = \Phi(2t) - \Phi(t) = \int_t^{2t} \phi(u)\mathrm{d}u > t\phi(2t) \gg t\phi(t)^2 = t^{-1}\Phi_1(t)^2,$$

as wanted. Now, if $V(t) \to \infty$ then the statement of the lemma follows from the above and the third estimate in Lemma 1. If $V_n \to \infty$ then $\Phi_{2n} - \Phi_n \to \infty$ since the mixed term in (3) is negative, and by the first estimate in Lemma 1 also $\Phi(2t) - \Phi(t) \to \infty$. $\square$

The rest of the argument for $K_n$ is as follows. From $\phi(t)^2 \ll \phi(2t)$ in the proof of Lemma 4 we get

$$\frac{\Phi(t+ct^{1/2}) - \Phi(t)}{(\Phi(2t) - \Phi(t))^{1/2}} = \frac{\int_t^{t+ct^{1/2}} \phi(u)\mathrm{d}u}{(\int_t^{2t} \phi(u)\mathrm{d}u)^{1/2}} \leq \frac{ct^{1/2}\phi(t)}{(t\phi(2t))^{1/2}} \to 0. \quad (13)$$

This implies

$$V(t+ct^{1/2})/V(t) \to 1 \quad (14)$$

provided $V(t) \to \infty$. Indeed,

$$\frac{V(t+ct^{1/2})}{V(t)} = \frac{\Phi(2t+2ct^{1/2}) - \Phi(t+ct^{1/2})}{\Phi(2t) - \Phi(t)} =$$

$$\frac{\Phi(2t+2ct^{1/2}) - \Phi(2t)}{\Phi(4t) - \Phi(2t)} \cdot \frac{\Phi(4t) - \Phi(2t)}{\Phi(2t) - \Phi(t)} - \frac{\Phi(t+ct^{1/2}) - \Phi(t)}{\Phi(2t) - \Phi(t)} + 1,$$

which tends to 1 by (13) and because $(\Phi(4t) - \Phi(2t))/(\Phi(2t) - \Phi(t)) < 2$ in consequence of $\Phi''(t) < 0$. From (13) and Lemma 1,

$$\mathbb{P}(K_{n+cn^{1/2}} - K_n > \epsilon V(n)^{1/2}) \leq \frac{\Phi_{n+cn^{1/2}} - \Phi_n}{\epsilon V(n)^{1/2}} = \frac{\Phi(n+cn^{1/2}) - \Phi(n)}{\epsilon V(n)^{1/2}} + o(1) \to 0.$$

In the very same way, and also using (14)

$$\mathbb{P}(K_n - K_{n-cn^{1/2}} > \epsilon V(n)^{1/2}) < \frac{\Phi(n) - \Phi(n-cn^{1/2})}{\epsilon V(n-cn^{1/2})^{1/2}} \left(\frac{V(n-cn^{1/2})}{V(n^{1/2})}\right)^{1/2} + o(1) \to 0.$$

Therefore $|K_{n\pm cn^{1/2}} - K_n|/V(n)^{1/2}$ converge to 0 in probability. Choosing $c$ sufficiently large we have for the Poisson process $n - cn^{1/2} < P(n) < n + cn^{1/2}$ and therefore $K_{n-cn^{1/2}} < K(n) < K_{n+cn^{1/2}}$ with probability larger $1 - \epsilon$. If follows that $(K_n - K(n))/V(n)^{1/2}$ converge to 0 in probability. The CLT for $K_n$ now follows from this and the CLT for $K(n)$. □

**Remark.** Hwang and Janson [25] have shown a more delicate local CLT for $K_n$ under $V(t) \to \infty$. If $\Phi_r(t)$ and $\text{Var}[K_r(t)]$ tend to $\infty$, then a CLT holds for $K_r(t)$, and one can naturally suspect that the same is valid for $K_{n,r}$ (this seems to have not been discussed in the literature). Mikhailov [31] proves a CLT for $K_{n,r}$ but assuming that $(p_j)$ vary with $n$ is a suitable way.

## 6. The variance

If $V(t)$ does not go to $\infty$ then $K(t)$ need not converge in distribution at all. Thus it is important to have criteria for $V(t) \to \infty$ and to understand other possible modes of bahaviour of the variance. For various $(p_j)$ the variance $V(t)$ and $V_n$ can go to $\infty$, converge to a finite limit, oscillate within a bounded range, or even oscillate between 0 to $\infty$.

In this section we sketch some recent results from [8]. The next lemma relates the variance with the mean number of singleton boxes.

**Lemma 5.** *There exists an increasing function $\tau(t)$ which satisfies $\tau(t)/t \in ]1, 2[$ for $t > 0$ and*

$$V(t) = \frac{t}{\tau(t)} \Phi_1(\tau(t)).$$

*Proof.* For $t > 0$ the equation $\Phi(2t) - \Phi(t) = t\Phi'(\tau)$ has a unique solution $\tau = \tau(t)$, which increases because $\Phi(t)$ is concave. □

It follows that $V(t) \to \infty$ is equivalent to $\Phi_1(t) \to \infty$. If $V(t)$ is bounded then $\Phi_1(t)$ is bounded and $V(n) - V_n \to 0$ by Lemma 1. We record some sufficient conditions for $V(t) \to \infty$.

**Proposition 6.** *Each of the following conditions implies $V(t) \to \infty$:*

(i) $\nabla \nu(x) \to \infty$ $(x \to 0)$, where $\nabla \nu(x) := \vec{\nu}(x/2) - \vec{\nu}(x) = \#\{j : x/2 \leq p_j < x\}$.
(ii) $\liminf_j p_{j+k}/p_j \geq 1/2$ for every $k = 1, 2, \ldots$
(iii) $\lim_{j \to \infty} p_j / \sum_{i \geq j} p_i = 0$.

*Proof.* Sufficiency of conditions (i) and (ii) is shown by rewriting (6) in the form

$$V(t) = \Phi(2t) - \Phi(t) = t \int_0^\infty e^{-tx} \nabla \nu(x) \mathrm{d}x. \tag{15}$$

For (iii) one exploits Lemma 5. □

All three conditions (i), (ii) and (iii) are satisfied (hence, $V(t)$ and $V_n$ converge to $\infty$) if $p_{j+1}/p_j \to 1$ $(j \to \infty)$.

We turn next to conditions for bounded variance or converging to a finite limit. The case of geometric frequencies gives a clue.

**Example 7.** Suppose $(p_j)$ is a geometric sequence $p_j = (1-q)q^j$ with ratio $0 < q < 1$. If $q = 1/2$ then $\nabla \nu(x) = 1$ for $0 < x < 1$ hence from (15) $V(t) \to 1$. More generally, for $q = 2^{-1/k}$ for some $k = 1, 2, \ldots$ we have $\nabla \nu(x) = k$ ($0 < x < 1$), hence $V(t) \to k$. For other values of $q$ the variance does not converge, rather it has an asymptotic expansion $V(t) = \log_{1/q} 2 + g(\log_{1/q} t) + o(1)$, where $g$ is a periodic function with mean 0 and a small amplitude [2].

**Proposition 8.** *A finite limit $v := \lim_{n \to \infty} V_n$ exists if and only if the frequencies satisfy*

$$\lim_{j \to \infty} \frac{p_{j+k}}{p_j} = \frac{1}{2}$$

*for some integer $k \geq 1$, in which case $v = k$.*

It can be shown [8] that the condition of proposition holds if and only if $(p_j)$ can be split in $k$ nonincreasing disjoint subsequence $(p_j^{(i)})$ $(i = 1, \ldots, k)$ such that for each of these subsequences $p_{j+1}^{(i)}/p_j^{(i)} \to 1/2$, hence the variance for sampling from $(cp_j^{(i)})$ approaches 1.

**Proposition 9.** $\limsup V(t) < \infty$ *if and only if there exists a positive integer $k$ such that for all $j = 1, 2, \ldots$*

$$\frac{p_{j+k}}{p_j} \leq \frac{1}{2}.$$

*Moreover, if for some $k \in \mathbb{N}$*

$$\limsup_{j \to \infty} \frac{p_{j+k}}{p_j} \leq \frac{1}{2} \tag{16}$$

*then $\limsup V(t) \leq k$, and this asymptotic bound is the best possible for frequencies satisfying (16) with given $k$.*

Many examples of irregular behaviour of $V(t)$ or $\Phi_r(t)$'s can be constructed using the following simple idea. Consider first a series of finitely many, say $m$, boxes with the same frequency $q$. In the Poisson scheme, the variance of the number of occupied boxes among these $m$ is $m(e^{-tq} - e^{-2tq})$ which is a unimodal function with the initial value 0, the maximum value $m/4$ assumed at $q^{-1} \log 2$, and exponential decay for larger $t$. Similar properties has the mean number of boxes occupied by $r$ balls, which is $mtqe^{-tq}$. Note that $m$ accounts for the maximum value, while varying $q$ amounts to just rescaling the time. Now, selecting $q_1 > q_2 > \ldots$ and taking $m_i$ frequencies equal $q_i$ ($i = 1, 2, \ldots$) we obtain a superposition of functions of the above type, thus creating oscillations with fairly arbitrary highs and lows. In the following examples we focus on the Poisson scheme, hence can ignore the normalization and only require $\sum_j p_j < \infty$.

**Example 10.** Choosing $q_i = q^i$ with some $0 < q < 1/2$ and $m_i = i$ we obtain a collection of frequencies $(p_j)$ for which $V(t) \to \infty$ ($t \to \infty$) but $\nabla \nu(x)$ oscillates between 0 and $\infty$. The example shows that condition (i) of Proposition 6 is not necessary for $V(t) \to \infty$.

**Example 11.** [28, p. 384] Choosing $q_i = 2^{-2^i}$ and $m_i = i$ we obtain $(p_j)$ for which $V(t_k) \to \infty$ along $t_k = 2^{2^k}$ but $V(t_k) \to 0$ along $t_k = 2^{2^k+k}$ ($k \to \infty$). Thus $V(t)$ oscillates between 0 and $\infty$, and the same applies to $\Phi_1(t)$.

**Example 12.** [8] Choosing $q_i = 2^{-2^{i+1}}$ and $m_i = 2^{2^i}$ we obtain $(p_j)$ for which $V(t) \to \infty$ and $\Phi_1(t) \to \infty$, while $\liminf \Phi_2(t) = 0$ and $\limsup \Phi_2(t) = \infty$. Thus we have here a curious pathology, when the mean number of singleton boxes goes to $\infty$, but the mean number of doubletons does not.

**Remark.** The last example disproves the assertion of [28, Lemma 1] that if the convergence radius of the series $\sum_j p_j u^j$ exceeds 1, then $\vec{\nu}(x)$ is slowly varying for $x \to 0$. Recall that slow variation means $\lim_{x \to 0} \vec{\nu}(cx)/\vec{\nu}(x) = 1$ for every $c > 0$. In the last example the convergence radius of $\sum_j p_j u^j$ equals 2, because $(p_j)^{1/j}$ oscillates between $1/4$ and $1/2$, as is readily checked. The multiplicity of each frequency $q_i = 2^{-2^i}$ is equal to the number of terms larger than this value ($k = 1, 2, \ldots$), which together with $c 2^{-2^k} < 2^{-2^{k-1}}$ (for fixed $c > 1$ and $k = k(c)$ large enough) implies that

$$\lim_{k \to \infty} \frac{\vec{\nu}(c 2^{-2^k})}{\vec{\nu}(2^{-2^k})} = \frac{1}{2} \quad (c > 1),$$

hence slow variation fails.

## 7. Regularly varying frequencies

We shall establish equivalent forms of regular variation in the occupancy problem, in particular we translate them in terms of the growth of mean values of $K_n$ and $K_{n,r}$'s.

Following [28] we say that the frequencies $(p_j)$ are *regularly varying* if

$$\vec{\nu}(x) \sim \ell(1/x) \, x^{-\alpha} \quad (x \downarrow 0) \tag{17}$$

for some $0 \leq \alpha \leq 1$ and a function $\ell$ slowly varying at $\infty$, i.e. satisfying $\ell(cy)/\ell(y) \to 1$ as $y \to \infty$, for every $c > 0$. The case $\alpha = 0$ corresponds to *slow variation*, while in the case $\alpha = 1$ we shall speak of *rapid variation*. In the case $\alpha = 1$ the summability of frequencies forces the function $\ell(1/x)$ to approach 0 (as $x \downarrow 0$) sufficiently fast.

Define for $r = 1, 2, \ldots$ the measures

$$\nu_r(\mathrm{d}x) := x^r \nu(\mathrm{d}x) = \sum p_j^r \delta_{p_j}(\mathrm{d}x),$$

so $\nu_r[0, x] = \sum_j p_j^r \mathbf{1}(p_j \leq x)$. The measure $\nu_1$ is the distribution of the frequency of the first discovered box, also called *the structural distribution* or *the law of the tagged fragment*, especially when the frequencies are random [6, 33].

**Proposition 13.** *The relation* (17) *with* $0 < \alpha < 1$ *is equivalent to*

$$\nu_1[0, x] \sim \frac{\alpha}{1-\alpha} x^{1-\alpha} \ell(1/x) \quad (x \downarrow 0), \tag{18}$$

*and for $r \geq 1$ it implies*

$$\nu_r[0, x] \sim \frac{\alpha}{r-\alpha} x^{r-\alpha} \ell(1/x) \quad (x \downarrow 0). \tag{19}$$

*Proof.* Integration by parts yields

$$\nu_1[0, x] = \int_0^x u\, \nu(\mathrm{d}u) = -x\, \vec{\nu}(x) + \int_0^x \vec{\nu}(u)\, \mathrm{d}u. \tag{20}$$

If (17) holds then, by Karamata's theorem [15, Theorem 1, Section 9 Ch. 8], the integral term is asymptotic to $(1-\alpha)^{-1} x^{1-\alpha} \ell(1/x)$, which readily yields (18). To show the converse implication we use

$$\vec{\nu}(x) = \int_x^\infty u^{-1}\, \nu_1(\mathrm{d}u) = c - x^{-1}\nu_1[0, x] + \int_x^\infty u^{-2}\nu_1[0, u]\, \mathrm{d}u, \tag{21}$$

evaluate the integral by Karamata's theorem and note that the constant term $c$ is dominated.

Similarly, applying Karamata's theorem to the integral term in

$$\nu_r[0, x] = \int_0^x u^r \nu(\mathrm{d}u) = -x^r \vec{\nu}(x) + r\int_0^x u^{r-1}\vec{\nu}(u)\mathrm{d}u,$$

we conclude (19). $\square$

The cases of slow and rapid variation need special treatment.

**Proposition 14.** *For $\ell$ slowly varying the relation*

$$\vec{\nu}(x) \sim x^{-1}\ell(1/x) \quad (x \downarrow 0) \tag{22}$$

*implies*

$$\nu_1[0, x] \sim \ell_1(1/x) \quad (x \downarrow 0), \tag{23}$$

with $\ell_1 \gg \ell$ another function of slow variation defined for $y > 1$ by

$$\ell_1(y) = \int_y^\infty u^{-1}\ell(u)\,\mathrm{d}u\,. \tag{24}$$

For $r > 1$ the relation (22) also implies

$$\nu_r[0, x] \sim \frac{1}{r-1} x^{r-1}\ell(1/x) \quad (x \downarrow 0). \tag{25}$$

In general, the relation (23) with some slowly varying $\ell_1$ only implies

$$\vec{\nu}(x) \ll x^{-1}\ell_1(1/x) \quad (x \downarrow 0), \tag{26}$$

and does not imply the regular variation of $\vec{\nu}(x)$; however if the regular variation holds then $\vec{\nu}(x)$ fulfills (22) with $\ell$ satisfying (24).

*Proof.* Suppose (22) is true. The integral (24) converges due to $\sum p_j = \int_0^\infty x\nu(\mathrm{d}x) < \infty$, and $\ell_1 \gg \ell$ since the integral of $u^{-1}$ diverges. The relation (23) follows from (20) by noting that the second term in the right side of (20) dominates the first. Asymptotics (25) follow by Karamata's theorem.

Conversely, (23) entails (26) by the virtue of (21). If $\vec{\nu}(x)$ is regularly varying then a direct argument shows that (22) and (24) must hold to match with (23). □

**Proposition 15.** *For $\ell_0$ slowly varying the relation*

$$\nu_1[0, x] \sim x\ell_0(1/x) \quad (x \downarrow 0) \tag{27}$$

*implies*

$$\vec{\nu}(x) \sim \ell(1/x) \quad (x \downarrow 0), \tag{28}$$

*with another slowly varying $\ell \gg \ell_0$ defined for $y > 1$ by*

$$\ell(y) = \int_1^y u^{-1}\ell_0(u)\mathrm{d}u\,, \tag{29}$$

*and also implies for $r > 1$*

$$\nu_r[0, x] \sim \frac{1}{r} x^r \ell_0(1/x) \quad (x \downarrow 0). \tag{30}$$

*In general, the relation (28) with some $\ell$ only implies*

$$\nu_1[0, x] \ll x\ell(1/x) \quad (x \downarrow 0)$$

*and does not imply the regular variation of $\nu_1[0, x]$; but if the regular variation holds, then (27) is satisfied with slowly varying $\ell_0$ related to $\ell$ via (29).*

*Proof.* The line of argument repeats the one in the previous proposition. Formula (30) is obtained by evaluating the terms in

$$\nu_r[0,x] = \int_0^x u^{r-1}\nu_1(\mathrm{d}u) = \nu_1[0,x]x^{r-1} - (r-1)\int_0^x \nu_1[0,u]u^{r-2}\mathrm{d}u.$$

□

The last proposition shows that (28) (i.e. (17) with $\alpha = 0$) is not strong enough to control $\nu_r$'s, but a slightly stronger assumption (27) is enough for that. The case of geometric frequencies demonstrates that (28) is indeed too weak.

**Example 16.** For $p_k = (1-q)q^{k-1}$, $k = 1, 2, \ldots$ ($0 < q < 1$) and $\nu = \sum_k \delta_{p_k}$ we have

$$\vec{\nu}(x) = \left\lceil \log_q \frac{x}{1-q} \right\rceil + 1 \sim |\log_q(1/x)| \quad (x \downarrow 0)$$

slowly varying, but

$$\nu_1[0,x] = q^{\lceil \log_q x/(1-q) \rceil}$$

is not regularly varying, since $x^{-1}\nu_1[0,x]$ oscillates between $(1-q)^{-1}$ and $q(1-q)^{-1}$.

We translate the above in terms of the mean values.

**Proposition 17.** *For $0 < \alpha < 1$, condition (17) is equivalent to each of the following two relations*

$$\Phi(t) \sim \Gamma(1-\alpha)t^\alpha \ell(t) \quad (t \to \infty), \tag{31}$$

*and*

$$\Phi_1(t) \sim \alpha\,\Gamma(1-\alpha)\,t^\alpha \ell(t) \quad (t \to \infty), \tag{32}$$

*and for $r \geq 1$ it implies*

$$\Phi_r(t) \sim \frac{\alpha\,\Gamma(r-\alpha)}{r!}\,t^\alpha \ell(t) \quad (t \to \infty). \tag{33}$$

*Proof.* Writing

$$\Phi(t) = t\int_0^\infty e^{-tx}\vec{\nu}(x)\mathrm{d}x,$$

we see that the equivalence of (31) and (17) follows by the Tauberian theorem for monotone densities [15, Theorem 4, Section 5, Ch. 13]. The equivalence of (32) and (18) follows from another form of the Tauberian theorem [15, Theorem 2, Section 5, Ch. 13]. In the same way, (33) follows from (19). □

In the case of rapid variation we have:

**Proposition 18.** *The relation* (22) *implies*

$$\Phi(t) \sim \Phi_1(t) \sim t\ell_1(t) \quad (t \to \infty), \tag{34}$$

*with $\ell_1$ as in* (24), *and for $r > 1$ it implies*

$$\Phi_r(t) \sim \frac{1}{r(r-1)}\, t\,\ell(t) \quad (t \to \infty). \tag{35}$$

*Also,* (23) *is equivalent to* $\Phi_1(t) \sim t\ell_1(t) \quad (t \to \infty)$.

*Proof.* Use Tauberian arguments and Proposition 14. □

In the case of slow variation we have:

**Proposition 19.** *The relation* (27) *implies* $\Phi(t) \sim \ell(t) \quad (t \to \infty)$, *and for $r \geq 1$*

$$\Phi_r(t) \sim \frac{1}{r}\ell_0(t) \quad (t \to \infty), \tag{36}$$

*where $\ell$ and $\ell_0$ are related as in* (29). *Also,* (27) *is equivalent to the $r = 1$ instance of* (36). *The relation* (28) *is equivalent to* $\Phi(t) \sim \ell(t) \quad (t \to \infty)$ *and entails* $\Phi_r(t) \ll \Phi(t) \quad (t \to \infty)$ *for every $r \geq 1$*.

*Proof.* Use Tauberian arguments and Proposition 15. □

We summarize the above relations for $\Phi(t)$ and $\Phi_r(t)$'s (the relations for $\Phi_n$ and $\Phi_{n,r}$'s are completely analogous):

**Corollary 20.** *In the case $0 < \alpha < 1$, for $r \geq 1$ all $\Phi_r(t)$ are of the same order of growth as $\Phi(t)$, and the ratios $\Phi_r(t)/\Phi(t)$ converge to $(-1)^r \binom{\alpha}{r}$ (these numbers comprise a probability distribution).*

*In the case of rapid variation $\Phi_r(t)$'s are of the same order for $r > 1$ but $\Phi(t) \sim \Phi_1(t) \gg \Phi_r(t)$ for $r > 1$, that is most of the occupied boxes are singleton.*

*In the case of slow variation under condition* (27) *we have $r^{-1}\Phi_r(t) \sim \Phi_1(t) \ll \Phi(t)$, meaning that all $\Phi_r(t)$'s are again of the same order but each of them is much smaller than $\Phi(t)$.*

**Remark.** The results about mean values are true also when $(p_j)$ are random, with $\nu$ understood as the intensity measure of the point process $\sum_j \delta_{p_j}$, and $\nu_1$ being the structural distribution. For instance, when $(p_j)$ are Poisson-Dirichlet($\theta$) frequencies we have $\nu_1(\mathrm{d}x) = \theta(1-x)^{\theta-1}$, so we are in the case of slow variation (27) with $\nu_1[0,x] \sim \theta x$, hence $\Phi_r(t)/\Phi_1(t) \to 1/r$ meaning that in the long run the mean number of balls in all $r$ton boxes is approximately the same, for each $r$. The last fact has been long known, especially for $\theta = 1$ in the context of random permutations which can be associated with a random sample from Poisson-Dirichlet(1), see [1].

We conclude as in [28]:

**Corollary 21.** *For $0 < \alpha < 1$ the condition of regular variation (17) is equivalent to any of the following conditions*

$$K_n \sim_{\text{a.s.}} \Gamma(1-\alpha) n^\alpha \ell(n) \quad (n \to \infty), \qquad K(t) \sim_{\text{a.s.}} \Gamma(1-\alpha) t^\alpha \ell(t) \quad (t \to \infty),$$

*and it implies for $r = 1, 2, \ldots$*

$$K_{n,r} \sim \frac{\alpha \, \Gamma(r-\alpha)}{r!} n^\alpha \ell(n) \quad (n \to \infty), \quad K_r(t) \sim \frac{\alpha \, \Gamma(r-\alpha)}{r!} t^\alpha \ell(t) \quad (t \to \infty).$$

*Proof.* Combine Proposition 2, the remarks after it, and Proposition 17. □

We stress that in this situation the strong laws for $K_r(t)$'s follow from the strong laws for increasing processes $\sum_{s>r} K_s(t)$ and the fact that $\sum_{s>r} \Phi_s(t) \asymp \Phi(t) \asymp t^\alpha \ell(t)$.

**Remark.** Characterizations of regular variation through behaviour of ratios of integrals with distinct kernels are known as Mercerian Tauberian theorems [7]. For instance, $\Phi_1(t)/\Phi(t) \to \alpha$ for some constant $0 < \alpha < 1$ implies (17).

When (17) holds with $0 < \alpha < 1$ the covariance matrix (7) becomes

$$\sigma_{r,s} = -\frac{\alpha \Gamma(r+s-\alpha)}{r!s!} 2^{\alpha-r-s}, \qquad r \neq s$$

$$\sigma_{r,r} = -\frac{\alpha \Gamma(2r-\alpha)}{r!r!} 2^{\alpha-2r} + \frac{\alpha \Gamma(r-\alpha)}{r!}$$

Using arguments similar to that in Section 5 Karlin [28, Theorem 5] showed that the array

$$\left( \frac{K_{n,r} - \mathbb{E}\left[K_{n,r}\right]}{(n^\alpha \ell(n))^{1/2}}, \ r = 1, 2, \ldots \right) \tag{37}$$

converges in distribution to a multivariate Gaussian array with zero mean and this covariance matrix. A similar result [28, theorem 5′] is valid also in the rapid variation case $\alpha = 1$, but $K_{n,1}$ requires a scaling different from the scaling of other $K_{n,r}$'s with $r > 1$, since $\text{Var}[K_{n,1}] \sim n\ell^*(n)$ is of larger order than $\text{Var}[K_{n,r}] \sim c_r n \ell(n)$ for $r > 1$.

## 8. Two uses of inversion

We recall further facts about regular variation. For function $h : \mathbb{R}_+ \to \mathbb{R}_+$ regularly varying at infinity with index $\alpha > 0$ there exists an asymptotic inverse function $g$ which is regularly varying with index $1/\alpha$ and satisfies $h(g(y)) \sim g(h(y)) \sim y$ $(y \to \infty)$. The function $g$ is unique up to the asymptotic equivalence, see [7, Proposition 1.5.12].

For $\ell$ a function of slow variation at infinity, its de Bruijn conjugate is another function of slow variation $\ell^\#$ satisfying

$$\ell(y)\ell^\#(y\ell(y)) \to 1, \quad \ell^\#(y)\ell(y\ell^\#(y)) \to 1 \quad (y \to \infty). \tag{38}$$

The de Bruijn conjugate exists and is uniquely determined up to asymptotic equivalence, see [7, Proposition 1.5.13]. For instance $(\log y)^\# \sim (\log y)^{-1}$. The following inversion formula is adapted from [7, Proposition 1.5.15].

**Lemma 22.** *Let h and g be asymptotic inverses of one another, then for $\alpha > 0$ and $\ell$ slowly varying*

$$h(y) \sim y^\alpha \ell(y) \quad (y \to \infty) \iff g(y) \sim y^{1/\alpha} (\ell^{1/\alpha}(y^{1/\alpha}))^\# \quad (y \to \infty). \tag{39}$$

Let $\ell^*$ be the reciprocal of the function appearing in the inversion formula (39):

$$\ell^*(y) := \frac{1}{\{\ell^{1/\alpha}(y^{1/\alpha})\}^\#}.$$

Keep in mind that $\ell^*$ depends on $\alpha$. This function allows one to formulate the property of regular variation directly in terms of individual frequencies.

**Proposition 23.** *The condition (17) with $0 < \alpha < 1$ is equivalent to*

$$p_j \sim \ell^*(j) \, j^{-1/\alpha} \quad (j \to \infty), \tag{40}$$

*and it implies*

$$\sum_{j>k} p_j \sim \frac{\alpha}{1-\alpha} \ell^*(k) \, k^{1-1/\alpha} \quad (k \to \infty). \tag{41}$$

*Proof.* Consider

$$h(y) = \vec{\nu}(1/y) \quad \text{and} \quad g(y) = \frac{1}{p_{\lceil y \rceil}}.$$

In view of

$$\vec{\nu}(1/y) = \max\left\{j : p_j \geq \frac{1}{y}\right\} = \sup\left\{z : \frac{1}{p_{\lceil z \rceil}} \leq y\right\},$$

the function $h$ is the generalized left-continuous inverse of $g$, and because $y \leq h(g(y)) < y + 1$, these functions are also asymptotic inverses of one another. By Lemma 22, the relation $h(y) \sim y^\alpha \ell(y)$ is equivalent to $g(y) \sim y^{1/\alpha}\{\ell^{1/\alpha}(y^{1/\alpha})\}^\#$, which is the same as (40). □

For fixed-$n$ allocation scheme define $N_k := \min\{n : K_n = k\}$ to be the times when new boxes are discovered. The analogous poissonized quantity is $T_k := \min\{t : K(t) = k\}$. By the definition

$$K_{N_k} = K_{T_k} = k \quad (k = 1, 2, \ldots).$$

Next proposition gives the relations inverse to that.

**Proposition 24.** *Under (17) with $0 < \alpha < 1$ we have*

$$N_k \sim_{\text{a.s.}} T_k \sim_{\text{a.s.}} \left(\frac{k}{\Gamma(1-\alpha)}\right)^{1/\alpha} \frac{1}{\ell^*(k)} \quad (k \to \infty).$$

*Moreover, the following asymptotic relations hold:*

$$n_k \sim \left(\frac{k}{\Gamma(1-\alpha)}\right)^{1/\alpha} \frac{1}{\ell^*(k)} \iff K_{n_k} \sim_{\text{a.s.}} K(n_k) \sim_{\text{a.s.}} k \quad (k \to \infty), \quad (42)$$

$$k_n \sim \Gamma(1-\alpha) n^\alpha \ell(n) \iff N_{k_n} \sim_{\text{a.s.}} T_{k_n} \sim_{\text{a.s.}} n \quad (n \to \infty). \quad (43)$$

*Proof.* Immediate from the uniqueness of the asymptotic inverse and de Bruijn conjugate, and Lemma 39. □

## 9. The cumulative frequency of empty boxes

Functionals $K_n$, $K_{n,r}$ are of the form

$$\sum_j \psi_j(X_{n,j}),$$

called 'separable statistics' by some authors [26, 17]. Here, $\psi_j$'s are some functions for which the sum is well defined. In this section we discuss one more instance of this kind,

$$S_n := \sum_j p_j \mathbf{1}(X_{n,j} = 0), \ S(t) := \sum_j p_j \mathbf{1}(X_j(t) = 0)$$

under the regular variation assumption (17) with $0 < \alpha < 1$. These have the meaning of the cumulative frequency of yet undiscovered boxes.

Defining $(\widetilde{p}_k)$ to be a random arrangement of frequencies in the order as the boxes are discovered, we can also write

$$S_n = \sum_{k > K_n} \widetilde{p}_k, \quad S(t) = \sum_{k > K(t)} \widetilde{p}_k, \quad R_k := \sum_{j > k} \widetilde{p}_j.$$

The sequence $(\widetilde{p}_k)$ is called a size-biased permutation of the frequencies $(p_j)$, see [33]. By (10) and (12),

$$\mathbb{E}[S_n] = n^{-1}\Phi(n, 1), \quad \mathbb{E}[S(t)] = t^{-1}\Phi_1(t).$$

These mean values control the geometric number of balls (respectively, the exponential time) needed to discover yet another box after some time. Clearly,

$$S_{N_k} = S(T_k) = R_k, \quad R_{K_n} = S_n, \quad R_{K(t)} = S(t). \quad (44)$$

**Proposition 25.** *Under assumption* (17) *with* $0 < \alpha < 1$

$$S_n \sim_{\text{a.s.}} S(n) \sim_{\text{a.s.}} \alpha \Gamma(1-\alpha) n^{\alpha-1} \ell(n) \quad (n \to \infty) \quad (45)$$

$$R_k \sim_{\text{a.s.}} \alpha \{\Gamma(1-\alpha)\}^{1/\alpha} k^{1-1/\alpha} \ell^*(k) \quad (k \to \infty). \quad (46)$$

*Proof.* From $S(t) = \sum_j \mathbf{1}(X_j(t) = 0) p_j$ using the same asymptotic evaluations as for (33) we compute the variance as

$$\text{Var}[S(t)] = \sum_j p_j^2 (e^{-p_j t} - e^{-2p_j t}) = \int_0^\infty x^2 (e^{-tx} - e^{-2tx}) \nu(\mathrm{d}x) \sim$$

$$\frac{\alpha \Gamma(2-\alpha)}{2} \left\{ t^{\alpha-2} \ell(t) - (2t)^{\alpha-2} \ell(2t) \right\} \sim \frac{\alpha \Gamma(2-\alpha)(1 - 2^{\alpha-2})}{2} t^{\alpha-2} \ell(t).$$

Thus we have

$$\mathbb{E}\left[S(t)\right] = t^{-1} \Phi_1(t) \sim c_1 t^{\alpha-1} \ell(t), \quad \text{Var}[S(t)] \sim c_2 t^{\alpha-2} \ell(t).$$

with some positive constants $c_1$, $c_2$. Using Chebyshev's inequality and the Borel-Cantelli lemma it is not hard to show that the convergence $S(t_m)/\mathbb{E}\left[S(t_m)\right] \to_{\text{a.s.}} 1$ is secured along a sequence $t_m \sim m^{2/\alpha}$. But then $S(t)/\mathbb{E}\left[S(t)\right] \to_{\text{a.s.}} 1$ follows for $t \to \infty$ by the usual sandwich argument, since $S(t)$ and $\mathbb{E}\left[S(t)\right]$ are nonincreasing and $t_m/t_{m+1} \to 1$ as $m \to \infty$. The rest of (45) follows by checking that with probability one $S(n(1+\epsilon)) < S_n < S(n(1-\epsilon))$ holds for all sufficiently large $n$.

To pass from (45) to (46) we first note that for $n_k$ as in (42)

$$S_{n_k} \sim_{\text{a.s.}} \alpha \{\Gamma(1-\alpha)\}^{1/\alpha} k^{1-1/\alpha} \ell^*(k) \quad (k \to \infty), \tag{47}$$

as is easily seen from

$$\{\ell^{1/\alpha}(y)\}^{\#} \ell^{1/\alpha} \left( y \{\ell^{1/\alpha}(y)\}^{\#} \right) \to 1,$$

which in turn is the second identity in (38). By Proposition 24 the bounds $n_{k(1-\epsilon)} < N_k < n_{k(1+\epsilon)}$ hold for sufficiently large $k$ with probability one. Finally, the first identity in (44) and monotonicity imply that for large $k$

$$S_{n_{k(1+\epsilon)}} < S_{N_k} = R_k < S_{n_{k(1-\epsilon)}},$$

and now (46) follows from (47) by letting $\epsilon \to 0$. $\square$

**Remark.** Comparing (46) with (41) we see that

$$R_k \Big/ \sum_{j>k} p_j \to_{\text{a.s.}} \{\Gamma(2-\alpha)\}^{1/\alpha} (1-\alpha)^{1-1/\alpha} > 1$$

which supports the intuition that the ranked arrangement of frequencies $(p_j)$ decays faster than the size-biased permutation $(\widetilde{p}_j)$. However, despite the fact that the tail sums of the sequences $(p_j)$ and $(\widetilde{p}_j)$ decay with the same order, the analogue of (40) for $(\widetilde{p}_j)$ does not hold. Indeed, by regular variation we have $p_i/p_{i+1} \to 1$ $(i \to \infty)$. Pick $\epsilon$ small, and let $p_i, p_{i+1}$ be the two largest frequencies smaller than $\epsilon$, and such that $p_i/p_{i+1}$ is close to 1. In the Poisson setup, let $a$ be the time needed to discover one of the two boxes with these frequencies $p_i$, $p_{i+1}$, and $b$ be the time to discover both of them. These can be expressed through independent exponential variables, showing that the ratio $a/b$ assumes values smaller than, say, $1/2$ with a probability bounded from zero. But this and the asymptotics of $K(t)$ readily imply that, with probability bounded away from $0$, the positions of $p_i$ and $p_{i+1}$ in the re-arranged sequence $(\widetilde{p}_j)$ are not asymptotic to one another, which could not happen if $(\widetilde{p}_j)$ were asymptotic to a regularly varying sequence.

## 10. Pure power laws

As has been mentioned in the Introduction, processes of coagulation and fragmentation of random masses [6, 33] are related to occupancy schemes where the frequencies $(p_j)$ are random. In many situations [22, 24, 32, 5] one encounters a relation

$$K_n \sim_{\text{a.s.}} D\, n^\alpha \qquad (n \to \infty) \tag{48}$$

for the number of blocks of partition induced by sampling from $(p_j)$, where $0 < \alpha < 1$ and $D$ is a strictly positive random variable (in this context, a measure of '$\alpha$-diversity' of the partition [33]). The asymptotics (48) say that, given $D$, $K_n$ is regularly varying with constant slow variation factor.

Conditioning on the frequencies and applying Proposition 2 we always have

$$K_n \sim_{\text{a.s.}} \mathbb{E}\left[K_n \mid (p_j)\right],$$

hence, by Proposition 17, (48) is equivalent to

$$\#\{j : p_j \geq x\} \sim_{\text{a.s.}} \frac{D}{\Gamma(1-\alpha)} x^{-\alpha} \qquad (x \to 0), \tag{49}$$

and by Proposition 23 it is also equivalent to

$$p_j \sim_{\text{a.s.}} \left(\frac{D}{\Gamma(1-\alpha)}\right)^{1/\alpha} j^{-\alpha} \qquad (j \to \infty). \tag{50}$$

Furthermore, any of these implies for $r = 1, 2, \ldots$

$$K_{n,r} \sim_{\text{a.s.}} \frac{\alpha(1-\alpha)\cdots(r-1-\alpha)}{r!}\, D n^\alpha \qquad (n \to \infty), \tag{51}$$

and for size-biased frequencies

$$\sum_{j>k} \widetilde{p}_j \sim_{\text{a.s.}} \alpha D^{1/\alpha} k^{1-1/\alpha} \qquad (k \to \infty). \tag{52}$$

The relations (48),(49),(50),(51),(52) appear in the literature under the folk name 'power laws'. Typically the starting point is (49), from which one arrives to the conclusion that $K_n/n^\alpha$ and $K_{n,r}/n^\alpha$ converge almost surely to multiples of the same random variable. But other direction is also useful: from the asymptotics of $K_n$, established by either the analysis of moments or some other method, one can make conclusions on the behaviour of small frequencies.

**Remark.** Interestingly, the distribution of $K_n$ does not converge to normal, although this is true for $K_n$ conditioned on $(p_j)$. An explanation for this phenomenon is that the randomness of $(p_j)$ dominates the variability due to random sampling. The CLT for $K_n$ does hold for sampling from Poisson-Dirichlet($\theta$) (in which case $K_n \sim_{a.s.} \theta \log n$), and has been also shown for some instances of random $(p_j)$ under more general assumptions of slow variation [4, 19, 23].

## 11. Strong laws for large parts

We collect here explicit distributional formulas and a few ways to describe the multivariate asymptotics of the 'large parts'.

For $X_n^\downarrow := (X_{n,j}^\downarrow, j = 1, 2, \ldots)$ the sequence of $X_{n,j}$'s arranged in nonincreasing order we have

$$\mathbb{P}(X_{n,j}^\downarrow = n_j, \ j = 1, \ldots, k) = \frac{n!}{\prod_{j=1}^k n_j!} \sum_{\substack{\text{distinct} \\ j_1, \ldots, j_k}} p_{j_1}^{n_1} \cdots p_{j_k}^{n_k}$$

$$(n_1 \geq \ldots \geq n_k > 0, \ n = \sum_{j=1}^k n_j),$$

where the sum expands over all $k$-tuples of distinct positive integers $j_1, \ldots, j_k$. Let $\widetilde{X}_n := (\widetilde{X}_{n,j}, j = 1, 2, \ldots)$ be the sequence which starts with positive terms of $X_n^\downarrow$ arranged in the order as the boxes are discovered by the balls, and ends with infinitely many zeroes. Similarly to the above,

$$\mathbb{P}(\widetilde{X}_{n,j}^\downarrow = n_j, \ j = 1, \ldots, k) = \frac{n!}{\prod_{j=1}^k (n_j - 1)!(n_j + n_{j+1} + \ldots + n_k)} \sum_{\substack{\text{distinct} \\ j_1, \ldots, j_k}} p_{j_1}^{n_1} \cdots p_{j_k}^{n_k},$$

where $n_1 > 0, \ldots, n_k > 0$, $n = n_1 + \ldots + n_k$. Let $(\widetilde{p}_j, j = 1, 2, \ldots)$ be the size-biased permutation of $(p_j)$. In particular, $\widetilde{X}_{n,1}$ is the number of balls out of the first $n$ which fall in the same box as the first ball, and $\widetilde{p}_1$ is the frequency of this box. The latter distribution conditionally given $(\widetilde{p}_j)$ is

$$\mathbb{P}(\widetilde{X}_{n,j}^\downarrow = n_j, \ j = 1, \ldots, k \,|\, (\widetilde{p}_j))$$
$$= \frac{n!}{\prod_{j=1}^k (n_j - 1)!(n_j + n_{j+1} + \ldots + n_k)} \prod_{j=1}^k (1 - \widetilde{p}_1 - \ldots - \widetilde{p}_{j-1})\widetilde{p}_j^{n_j - 1}.$$

To proceed to the asymptotics, recall that the succession of hits in box $j$ undergoes a Bernoulli process with success probability $p_j$, hence by the classical law of large numbers $n^{-1} X_{n,j} \to_{a.s.} p_j$ as $n \to \infty$. The following multivariate extensions of this result involve various arrangements of the boxes, and feature the behaviour of 'large' parts which grow linearly with $n$. These results are of fundamental importance in Kingman's theory of exchangeable partitions [6, 33].

**Proposition 26.** *As $n \to \infty$ we have*

$$n^{-1}X_n \to_{\text{a.s.}} (p_1, p_2, \ldots), \tag{53}$$
$$n^{-1}X_n^{\downarrow} \to_{\text{a.s.}} (p_1, p_2, \ldots), \tag{54}$$
$$n^{-1}\widetilde{X}_n \to_{\text{a.s.}} (\widetilde{p}_1, \widetilde{p}_2, \ldots), \tag{55}$$

*where the convergence is understood in the product topology.*

*Proof.* The first relation amounts to the marginal convergence. The second relation follows from $\sum_j n^{-1}X_{n,j} = 1$ and the fact that ranking is a continuous mapping on the infinite-dimensional simplex $\{(x_1, x_2, \ldots) : x_j \geq 0, \sum_j x_j = 1\}$. The third relation is a known consequence of the second [33]. □

Many questions on random allocations also involve the ordering of the boxes. One can be interested in the last occupied box [13] or the first empty [11] etc. This theme lies outside the scope of this paper, and we only mention one generalization of (53) that appears in the theory of exchangeable *ordered* partitions [18, 21].

Let $\triangleleft$ be a strict total order on positive integers, thought of as some 'arrangement' of boxes. With each box $j$ we associate an open interval $]a_j, b_j[$ of length $p_j$, with endpoints $a_j = \sum_{\{i: i \triangleleft j\}} p_i$ and $b_j = a_j + p_j$. Let $\mathcal{C} = [0, 1] \setminus (\cup_j ]a_j, b_j[)$ be the complementary closed set. For each $n$ define a random finite set $\mathcal{C}_n$ comprised of 0 and of all distinct elements of the sequence $(\sum_{\{i: i \triangleleft j\}} n^{-1}X_{n,i}, j = 1, 2, \ldots)$.

**Proposition 27.** *For $d_H$ the Hausdorff distance on the space of closed subsets of $[0, 1]$,*

$$d_H(\mathcal{C}_n, \mathcal{C}) \to_{\text{a.s.}} 0 \quad (n \to \infty).$$

*Proof.* This follows from (54) and the easily established convergence $\sum_{\{i: i \triangleleft j\}} n^{-1}X_{n,i} \to_{\text{a.s.}} a_j$. □

For instance, define the order $\triangleleft$ by choosing a sequence of distinct reals $(y_j, j = 1, 2, \ldots)$ and setting $i \triangleleft j$ iff $y_i < y_j$. Consider the discrete distribution $\mu(\mathrm{d}y) = \sum_j p_j \delta_{y_j}(\mathrm{d}y)$ which places mass $p_j$ at point $y_j$. For an independent sample of $n$ elements from $\mu$, the finite set $\mathcal{C}_n$ encodes the nonzero counts of repeated values in the sample, in the natural (increasing) order of the values. The set $\mathcal{C}$ can be identified with the quantile transform of $\mu$ [18]. For example, $\triangleleft$ is the standard order on integers for $y_j = j$ ($j = 1, 2, \ldots$), in which case $\mathcal{C} = \{0, p_1, p_1 + p_2, p_1 + p_2 + p_3, \ldots, 1\}$ and Proposition 27 amounts to (54). A more sophisticated example appears when $\{y_j\}$ is the set of all rational numbers, in which case $\mathcal{C}$ is a Cantor set.

**Acknowledgement** We are indebted to a referee for constructive comments.